\documentclass[twoside,11pt]{amsart}
\usepackage{amsmath,amssymb,amsthm}

\begin{document}

\setlength{\textwidth}{145mm} \setlength{\textheight}{203mm}
\setlength{\parindent}{0mm} \setlength{\parskip}{2pt plus 2pt}

\frenchspacing



\numberwithin{equation}{section}
\newtheorem{thm}{Theorem}[section]
\newtheorem{lem}[thm]{Lemma}
\newtheorem{prop}[thm]{Proposition}
\newtheorem{cor}[thm]{Corollary}
\newtheorem{probl}[thm]{Problem}

\newtheorem{defn}{Definition}[section]
\newtheorem{rem}{Remark}[section]
\newtheorem{exa}{Example}

\newcommand{\be}[1]{\begin{equation}\label{#1}}
\newcommand{\ee}{\end{equation}}


\newcommand{\X}{\mathfrak{X}}
\newcommand{\B}{\mathcal{B}}
\newcommand{\s}{\mathfrak{S}}
\newcommand{\g}{\mathfrak{g}}
\newcommand{\W}{\mathcal{W}}
\newcommand{\T}{\mathcal{T}}
\newcommand{\Lgr}{\mathrm{L}}
\newcommand{\dd}{\mathrm{d}}
\newcommand{\n}{\nabla}
\newcommand{\nn}{\nabla'}
\newcommand{\lm}{\lambda}
\newcommand{\ta}{\theta}
\newcommand{\pd}{\partial}
\newcommand{\LL}{\mathcal{L}}

\newcommand{\diag}{\mathrm{diag}}
\newcommand{\End}{\mathrm{End}}
\newcommand{\im}{\mathrm{Im}}
\newcommand{\id}{\mathrm{id}}

\newcommand{\ie}{i.e. }
\newfont{\w}{msbm9 scaled\magstep1}
\def\R{\mbox{\w R}}
\newcommand{\norm}[1]{\left\Vert#1\right\Vert ^2}
\newcommand{\nnorm}[1]{\left\Vert#1\right\Vert ^{*2}}
\newcommand{\nN}{\norm{N}}
\newcommand{\nP}{\norm{\nabla P}}
\newcommand{\nnP}{\nnorm{\nabla P}}
\newcommand{\tr}{{\rm tr}}

\newcommand{\thmref}[1]{Theorem~\ref{#1}}
\newcommand{\propref}[1]{Proposition~\ref{#1}}
\newcommand{\corref}[1]{Corollary~\ref{#1}}
\newcommand{\secref}[1]{\S\ref{#1}}
\newcommand{\lemref}[1]{Lemma~\ref{#1}}
\newcommand{\dfnref}[1]{Definition~\ref{#1}}

\frenchspacing


\title{Natural connections on Riemannian product manifolds}

\author{Dobrinka Gribacheva}

\maketitle

{\small
\textbf{Abstract} \\
A Riemannian almost product manifold with integrable almost
product structure is called a Riemannian product manifold. In the
present paper the natural connections on such manifolds are
studied, i.e. the linear connections preserving the almost product
structure and the Riemannian metric.
}


\section{Introduction}

The systematic development of the theory of Riemannian almost
product manifolds was started by K. Yano in \cite{Ya}, where basic
facts of the differential geometry of these manifolds are given. A
Riemannian almost product manifold $(M,P,g)$ is a differentiable
manifold $M$ for which almost product structure $P$ is compatible
with the Riemannian metric $g$ such that a isometry is induced in
any tangent space of $M$.

The geometry of a Riemannian almost product manifold $(M,P,g)$ is
a geometry of both structures $g$ and $P$. There are important in
this geometry the linear connections with respect to which the
parallel transport determine an isomorphism of the tangent spaces
with the structures $g$ and $P$. This is valid if and only if the
structures $g$ and $P$ are parallel with respect to such a
connection. In the general case on a Riemannian almost product
manifold there exist a countless number of linear connections
regarding which $g$ and $P$ are parallel. Such connections are
called natural in \cite{Mi}.

In the present work we consider the natural connections on the
Riemannian product manifolds $(M,P,g)$, i.e. on the Riemannian
almost product manifolds $(M,P,g)$ with an integrable structure
$P$. In our investigations we suppose the condition $\tr P=0$,
which implies that $\dim M$ is an even number.

In \cite{Nav} A.~M.~Naveira gave a classification of Riemannian
almost product manifolds with respect to the covariant
differentiation $\n P$, where $\n$ is the Levi-Civita connection
of $g$. Having in mind the results in \cite{Nav}, M.~Staikova and
K.~Gribachev gave in \cite{Sta-Gri} a classification of the
Riemannian almost product manifolds $(M,P,g)$ with $\tr P=0$.

In Section 2 we give some necessary facts about the Riemannian
almost product manifolds. We recall the classification of
Staikova-Gribachev for these manifolds which is made regarding the
tensor $F$ determined by $F(x,y,z)=g((\n_x P)y,z)$. The basic
classes are $\W_1$, $\W_2$ and $\W_3$. The class of the Riemannian
 product manifolds is $\W_1\oplus\W_2$.

In Section 3 we recall a decomposition of the space of the torsion
tensors on a Riemannian almost product manifold to invariant
orthogonal subspaces $\T_i$ $(i=1,2,3,4)$ given in \cite{Mi}. We
establish some properties of the torsion of a natural connection
on a manifold $(M,P,g)\in \W_1\oplus\W_2$ in terms of the
mentioned decomposition.

In Section 4 we establish that the unique natural connection on
$(M,P,g)\in \W_1\oplus\W_2$ with torsion $T$, which can be
expressed by the components of $F$, is the canonical connection.
We prove that this is the unique natural connection for which
$T\in\T_1$.

In Section 5 we consider the natural connections for which the
torsion $T$ can be expressed by the components of the tensor
$g\otimes\theta$, where $\theta$ is the Lee 1-form associated with
$F$. Such connections exist only on a manifold $(M,P,g)\in \W_1$
and their torsions belong to a 2-parametric family. When the
natural connection does not coincide with the canonical connection
then we have $T\in\T_1\oplus\T_4$, $T\notin\T_1$ and
$T\notin\T_4$.

\section{Preliminaries}

Let $(M,P,g)$ be a \emph{Riemannian almost product manifold},
\ie{} a differentiable manifold $M$ with a tensor field $P$ of
type $(1,1)$ and a Riemannian metric $g$ such that
\begin{equation*}\label{Pg}
    P^2x=x,\quad g(Px,Py)=g(x,y)
\end{equation*}
for arbitrary $x$, $y$ of the algebra $\X(M)$ of the smooth vector
fields on $M$. Obviously $g(Px,y)=g(x,Py)$.

Further $x,y,z,w$ will stand for arbitrary elements of $\X(M)$ or
vectors in the tangent space $T_pM$ at $p\in M$.

In this work we consider Riemannian almost product manifolds with
$\tr{P}=0$. In this case $(M,P,g)$ is an even-dimensional
manifold. We denote $\dim{M}=2n$.

The classification in \cite{Sta-Gri} of Riemannian almost product
manifolds is made with respect to the tensor field $F$ of type
(0,3), defined by
\begin{equation*}\label{2}
F(x,y,z)=g\left(\left(\nabla_x P\right)y,z\right),
\end{equation*}
where $\nabla$ is the Levi-Civita connection of $g$. The tensor
$F$ has the following properties:
\begin{equation}\label{3}
\begin{array}{l}
    F(x,y,z)=F(x,z,y)=-F(x,Py,Pz),\\[6pt] F(x,y,Pz)=-F(x,Py,z).
\end{array}
\end{equation}

\emph{The associated 1-form} $\theta$ for $F$ is determined by the
equality
\[
\theta(x)=g^{ij}F(e_i,e_j,x),
\]
where $g^{ij}$ are the components of the inverse matrix of $g$
with respect to the basis $\{e_i\}$ of $T_pM$.

The basic classes of the classification in \cite{Sta-Gri} are
$\W_1$, $\W_2$ and $\W_3$. Their intersection is the class $\W_0$
of the \emph{Riemannian $P$-manifolds}, determined by the
condition $F(x,y,z)=0$ or equivalently $\n P=0$. In the
classification there are include the classes $\W_1\oplus\W_2$,
$\W_1\oplus\W_3$, $\W_2\oplus\W_3$ and the class
$\W_1\oplus\W_2\oplus\W_3$ of all Riemannian almost product
manifolds.

In the present work we consider only the Riemannian almost product
manifolds $(M,P,g)$ with integrable almost product structure $P$,
i.e. the manifolds with zero Nijenhuis tensor $N$ determined by
\begin{equation*}\label{N}
    N(x,y)=\left(\nabla_x P\right)Py-\left(\nabla_y P\right)Px
    +\left(\nabla_{Px} P\right)y-\left(\nabla_{Py} P\right)x.
\end{equation*}
These manifolds is called \emph{Riemannian product manifolds} and
they form the class $\W_1\oplus\W_2$. The characteristic
conditions for the classes $\W_1$, $\W_2$ and $\W_1\oplus\W_2$ are
the following
\[
\begin{split}
\W_{1}:\qquad    & F(x,y,z)=\frac{1}{2n}\big\{ g(x,y)\ta (z)+g(x,z)\ta (y)-g(x,Py)\ta (Pz) \big.\\
& \phantom{F(x,y,z)=\frac{1}{2n}\big\{g(x,y)\ta (z)+g(x,z)\ta (y)
\big.} -g(x,Pz)\ta (Py)\big\};
\\
\W_{2}:\qquad      &
F(x,y,Pz)+F(y,z,Px)+F(z,x,Py)=0,\hspace{0.1in} \ta=0;
\\
\W_{1}\oplus \W_{2}:\qquad      &
F(x,y,Pz)+F(y,z,Px)+F(z,x,Py)=0.
\end{split}
\]


\section{Natural connections on Riemannian product manifolds}

The linear connections in our investigations have a torsion.

Let $\nn$ be a linear connection with a tensor $Q$ of the
transformation $\n \rightarrow\nn$ and a torsion $T$, \ie{}
\[
\nn_x y=\n_x y+Q(x,y),\quad T(x,y)=\nn_x y-\nn_y x-[x,y].
\]
The corresponding (0,3)-tensors are defined by
\begin{equation*}\label{3.1}
    Q(x,y,z)=g(Q(x,y),z), \quad T(x,y,z)=g(T(x,y),z).
\end{equation*}
The symmetry of the Levi-Civita connection implies
\begin{equation}\label{3.2}
    T(x,y)=Q(x,y)-Q(y,x),
\end{equation}
\[
    T(x,y)=-T(y,x).
\]

A partial decomposition of the space $\mathcal{T}$ of the torsion
tensors $T$ of type (0,3) is valid on a Riemannian  almost product
manifold $(M,P,g)$:
$\mathcal{T}=\mathcal{T}_1\oplus\mathcal{T}_2\oplus\mathcal{T}_3\oplus\mathcal{T}_4$,
where $\mathcal{T}_i$ $(i=1,2,3,4)$ are invariant orthogonal
subspaces \cite{Mi}. For the projection operators $p_i$ of
$\mathcal{T}$ in $\mathcal{T}_i$ is established:
\begin{equation*}
  \begin{split}
&
    p_1(x,y,z)=\frac{1}{8}\bigl\{2T(x,y,z)-T(y,z,x)-T(z,x,y)-T(Pz,x,Py)\bigr.\\[4pt]
& \phantom{p_1(x,y,z)=\frac{1}{8}}
    +T(Py,z,Px)+T(z,Px,Py)-2T(Px,Py,z)\\[4pt]
& \phantom{p_1(x,y,z)=\frac{1}{8}}
    +T(Py,Pz,x)+T(Pz,Px,y)-T(y,Pz,Px)\bigr\},\\[4pt]
&
    p_2(x,y,z)=\frac{1}{8}\bigl\{2T(x,y,z)+T(y,z,x)+T(z,x,y)+T(Pz,x,Py)\bigr.\\[4pt]
&  \phantom{p_2(x,y,z)=\frac{1}{8}}
    -T(Py,z,Px)-T(z,Px,Py)-2T(Px,Py,z)\\[4pt]
& \phantom{p_2(x,y,z)=\frac{1}{8}}
    -T(Py,Pz,x)-T(Pz,Px,y)+T(y,Pz,Px)\bigr\},\\[4pt]
  &
    p_3(x,y,z)=\frac{1}{4}\bigl\{T(x,y,z)+T(Px,Py,z)-T(Px,y,Pz)-T(x,Py,Pz)\bigr\},\\[4pt]
  &
    p_4(x,y,z)=\frac{1}{4}\bigl\{T(x,y,z)+T(Px,Py,z)+T(Px,y,Pz)+T(x,Py,Pz)\bigr\}.
  \end{split}
\end{equation*}

\begin{defn}[\cite{Mi}]\label{defn-3.1'}
A linear connection $\nn$ on a Riemannian almost product manifold
$(M,P,g)$ is called a \emph{natural connection} if $\nn P=\nn
g=0$.
\end{defn}
If $\nn$ is a linear connection with a tensor $Q$ of the
transformation $\n \rightarrow\nn$ on a Riemannian almost product
manifold, then it is  a natural connection if and only if the
following conditions are valid \cite{Mi}:
\begin{equation}\label{3.5}
    F(x,y,z)=Q(x,y,Pz)-Q(x,Py,z),
\end{equation}
\begin{equation}\label{3.6}
    Q(x,y,z)=-Q(x,z,y).
\end{equation}

Let $\Phi$ be the (0,3)-tensor determined by
\[
    \Phi(x,y,z)=g\left(\widetilde{\nabla}_x y - \n_x y,z \right),
\]
where $\widetilde{\nabla}$ is the Levi-Civita connection of the
\emph{associated metric} $\tilde{g}$ determined by
$\tilde{g}(x,y)=g(x,Py)$.

\begin{thm}[\cite{Mi}]\label{t-3.1}
A linear connection with the torsion  $T$ on a Riemannian almost
product manifold $(M,P,g)$ is natural if and only if
\begin{equation*}\label{3.9}
  \begin{array}{l}
    4p_1(x,y,z)=-\Phi(x,y,z)+\Phi(y,z,x)-\Phi(x,Py,Pz)\\[4pt]
    \phantom{4p_1(x,y,z)=}-\Phi(y,Pz,Px)+2\Phi(z,Px,Py),
  \end{array}
\end{equation*}
\begin{equation*}\label{3.8}
    4p_3(x,y,z)=-g(N(x,y),z)=-2\left\{\Phi(z,Px,Py)+\Phi(z,x,y)\right\}.
\end{equation*}
\end{thm}


Let $(M,P,g)$ be a Riemannian almost product manifold, i.e.
$(M,P,g)\in\W_1\oplus\W_2$. For such a manifold is valid the
following equality (\cite{Sta-Gri})
\begin{equation}\label{10}
    \Phi(x,y,z)=\frac{1}{2}\bigl\{F(y,x,Pz)-F(Py,x,z)\bigr\}.
\end{equation}

By virtue of \eqref{10}, the characteristic condition for the
class $\W_1\oplus\W_2$ and \thmref{t-3.1} we obtain the following
\begin{thm}\label{t-3.2}
A linear connection with the torsion  $T$ on a Riemannian product
manifold $(M,P,g)$ is natural if and only if
\begin{equation}\label{11}
    p_1(x,y,z)=\frac{1}{2}F(z,y,Px),\qquad p_3(x,y,z)=0.
\end{equation}
\qed
\end{thm}

According to \thmref{t-3.2} and the conditions for the projection
operators $p_2$ and $p_4$, we get the following
\begin{cor}\label{c-3.3}
For the torsion $T$ of a natural connection on a Riemannian
product manifold $(M,P,g)$ are valid the following equalities
\begin{gather}
    p_2(x,y,z)=\frac{1}{2}\left\{T(x,y,z)-T(Px,Py,z)+F(z,x,Py)\right\},\label{12}\\[4pt]
    p_4(x,y,z)=\frac{1}{2}\left\{T(x,y,z)+T(Px,Py,z)\right\}.\label{13}\nonumber
\end{gather}
\qed
\end{cor}

Further, we suppose that the considered Riemannian product
manifold \allowbreak $(M,P,g)$ is not a Riemannian $P$-manifold,
i.e. $F$ is not a zero tensor.

According to \thmref{t-3.2} for the torsion  $T$ of a natural
connection on a Riemannian product manifold $(M,P,g)$, we have
$T\in\T_1\oplus\T_2\oplus\T_4$. If we suppose that $T=p_2+p_4$
then, having in mind \corref{c-3.3}, we obtain $F=0$, which is a
contradiction. Therefore $T\notin\T_2\oplus\T_4$. Then we have to
consider the cases:
\begin{enumerate}
    \item[A)] $T\in\T_1$;
    \item[B)] $T\in\T_1\oplus\T_4$, $T\notin\T_1$, $T\notin\T_4$;
    \item[C)] $T\in\T_1\oplus\T_2$, $T\notin\T_1$, $T\notin\T_2$.
\end{enumerate}


\section{Case A}

In \cite{Mi} a natural connection on a Riemannian almost product
manifold \allowbreak $(M,P,g)$ is called \emph{canonical} if for
its torsion the following equality is valid
\begin{equation*}\label{14}
    T(x,y,z)+T(y,z,x)+T(Px,y,Pz)+T(y,Pz,Px)=0.
\end{equation*}
This connection is an analogue of the Hermitian connection on the
Hermitian manifolds (\cite{Li-1}). In \cite{Mi} it is proved that
on any Riemannian product manifold $(M,P,g)$ there exist a unique
canonical connection and for the torsion $T$ of this connection
the condition $T\in\T_1\oplus\T_3$ is valid. Then, according to
\thmref{t-3.2} it is valid the following
\begin{thm}\label{t-4.1}
The case A for the torsion $T$ of a natural connection on a
Riemannian product manifold $(M,P,g)$ is valid if and only if this
connection is the canonical one. In this case the following
equality is satisfied
\begin{equation}\label{15}
    T(x,y,z)=\frac{1}{2}F(z,y,Px).
\end{equation}
\qed
\end{thm}

Let $T$ is the torsion $T$ of a natural connection on a Riemannian
product manifold $(M,P,g)$. Having in mind the characteristic
condition for the class $\W_1\oplus\W_2$ and conditions \eqref{3},
we obtain the following expression of $T$ by the independent
components of $F$:
\begin{equation}\label{16}
\begin{split}
    T(x,y,z)&=\lm_1F(x,y,z)+\lm_2F(y,z,x)+\lm_3F(Px,y,z)\\[4pt]
            &+\lm_4F(Py,z,x)+\lm_5F(x,y,Pz)+\lm_6F(y,z,Px)\\[4pt]
            &+\lm_7F(Px,Py,z)+\lm_8F(Py,Pz,x),
\end{split}
\end{equation}
where $\lm_i\in\R$ $(i=1,2,\dots,8)$. From \eqref{16}, using
\eqref{3}, \eqref{3.2} and \eqref{3.6}, we get \eqref{15} and
therefore $T$ is the torsion of the canonical connection on
$(M,P,g)$.

Hence we establish that it is valid the following
\begin{prop}\label{p-4.2}
The canonical connection is the unique natural connection on a
Riemannian product manifold $(M,P,g)$, which torsion can be
expressed by the tensor $F$. \qed
\end{prop}

The canonical connection on a Riemannian product manifold
$(M,P,g)\in\W_1$ is studied in \cite{Sta-Gri}. Having in mind the
characteristic condition for the class $\W_1$ and condition
\eqref{15} for the torsion of the canonical connection on
$(M,P,g)\in\W_1\oplus\W_2$, we obtain the following
\begin{prop}\label{p-4.3}
For the torsion $T$ of the canonical connection on a Riemannian
product manifold $(M,P,g)\in\W_1$ the following equality is
valid
\begin{equation}\label{16'}
\begin{split}
    T(x,y,z)=\frac{1}{4n}\bigl\{&g(y,z)\ta(Px)-g(y,Pz)\ta(x)\\[4pt]
                          -&g(x,z)\ta(Py)+g(x,Pz)\ta(y)\bigr\}.
\end{split}
\end{equation}

 \qed
\end{prop}


\section{Case B and Case C}

Having in mind the latter two propositions, in the present section
for the cases B and C  we consider the existence of a natural
connection with torsion $T$ expressed by the components of the
tensor $g\otimes\theta$, i.e.
\begin{equation}\label{17}
\begin{split}
    T(x,y,z)&=\lm_1g(x,y)\ta(z)+\lm_2g(y,z)\ta(x)+\lm_3g(z,x)\ta(y)\\[4pt]
            &+\lm_4g(x,y)\ta(Pz)+\lm_5g(y,z)\ta(Px)+\lm_6g(z,x)\ta(Py)\\[4pt]
            &+\lm_7g(x,Py)\ta(z)+\lm_8g(y,Pz)\ta(x)+\lm_9g(z,Px)\ta(y)\\[4pt]
            &+\lm_{10}g(x,Py)\ta(Pz)+\lm_{11}g(y,Pz)\ta(Px)\\[4pt]
            &+\lm_{12}g(z,Px)\ta(Py).
\end{split}
\end{equation}

From \eqref{17}, using \eqref{3.2}, \eqref{3.5} and \eqref{3.6},
we obtain
\begin{equation*}\label{18}
\begin{split}
    &F(x,y,z)=\\[4pt]
    &\phantom{+}(\lm_1-\lm_{10})\left\{g(x,y)\ta(Pz)-g(x,Pz)\ta(y)-g(x,Py)\ta(z)+g(x,z)\ta(Py)\right\}\\[4pt]
            &+(\lm_4-\lm_7)\left\{g(x,y)\ta(z)-g(x,Pz)\ta(Py)-g(x,Py)\ta(Pz)+g(x,z)\ta(y)\right\}.
\end{split}
\end{equation*}

Since the tensor $F$ is expressed by the tensor $g\otimes\ta$ only
for the class $\W_1$, then the comparison of the latter equality
with the characteristic condition of the mentioned class implies
\[
    \lm_1=\lm_{10},\quad \lm_4-\lm_7=\frac{1}{2n},\quad
    \lm_2=\lm_5=\lm_6=\lm_8=\lm_9=\lm_{11}=\lm_{12}=0.
\]
Then from \eqref{17}, using the denotations $\lm=\lm_1=\lm_{10}$
and $\mu=\lm_7$, we obtain the following
\begin{thm}\label{t-5.1}
Let the torsion $T$ of a natural connection on a Riemannian
product manifold $(M,P,g)\in\W_1$ is expressed by
$g\otimes\theta$. Then Case B or Case C is valid for $T$ if and
only if $(M,P,g)\in\W_1$. In this case $T$ has the following
representation
\begin{equation}\label{19}
\begin{split}
    &T(x,y,z)=\\[4pt]
    &\phantom{+}\lm\left\{g(y,z)\ta(x)-g(x,z)\ta(y)+g(y,Pz)\ta(Px)-g(x,Pz)\ta(Py)\right\}\\[4pt]
            &+\mu\left\{g(y,Pz)\ta(x)-g(x,Pz)\ta(y)+g(y,z)\ta(Px)-g(x,z)\ta(Py)\right\}\\[4pt]
            &+\frac{1}{2n}\left\{g(y,z)\ta(Px)-g(x,z)\ta(Py)\right\},\qquad \lm,\mu\in\R.
\end{split}
\end{equation}
\qed
\end{thm}

Let us consider Case B, i.e. $T=p_1+p_4$. Then, according to
\eqref{11} and \eqref{12}, we have that this case is valid if and
inly if
\begin{equation}\label{20}
    F(z,x,Py)=T(Px,Py,z)-T(x,y,z).
\end{equation}
We verify directly that condition \eqref{20} is satisfied for any
torsion $T$ determined by \eqref{19}, i.e. for arbitrary $\lm$ and
$\mu$. Let us remark that for $\lm=0$ and $\mu=-\frac{1}{4n}$ from
\eqref{19} we get condition \eqref{16'} for the torsion of the
canonical connection on $(M,P,g)\in\W_1$, i.e. Case A for the
class $\W_1$. Therefore, it is valid the following
\begin{thm}\label{t-5.2}
Let the torsion $T$ of a natural connection on a Riemannian
product manifold $(M,P,g)\in\W_1$ is expressed by
$g\otimes\theta$. Then Case B is valid if and only if $T$ is
determined by \eqref{19} for $(\lm,\mu)\neq(0,-\frac{1}{4n})$.
\qed
\end{thm}

Let us consider Case C. According to \thmref{t-5.1}, the torsion
$T$ is determined by \eqref{19}. Then, having in mind
\thmref{t-5.2}, we establish that  $T$ satisfies the conditions of
Case B. Therefore, we obtain the following
\begin{prop}\label{p-5.3}
If the torsion $T$ of a natural connection on a Riemannian product
manifold $(M,P,g)\in\W_1$ is expressed by $g\otimes\theta$, then
Case C does not exist. \qed
\end{prop}


\section{Conclusion}

The canonical connection is the unique natural connection on any
Riemannian product manifold $(M,P,g)$, which torsion can be
expressed by $F$. This is the unique natural connection with
torsion in Case A.

If a natural connection on $(M,P,g)\in\W_1$ has a torsion $T$
expressed by $g\otimes\theta$, then $T$ belongs to a 2-parametric
family determined by \eqref{19}. Case A and Case B for $T$ are
valid  when $(\lm,\mu)=(0,-\frac{1}{4n})$ and
$(\lm,\mu)\neq(0,-\frac{1}{4n})$, respectively.

Since $\ta=0$ for the class $\W_2$, then there do not exist any
natural connection on $(M,P,g)\in\W_2$ with torsion expressed by
$g\otimes\theta$. Only Case A is valid on such a manifold.


\bigskip

\textit{Dobrinka Gribacheva\\
Department of Geometry\\
Faculty of Mathematics and Informatics
\\
Paisii Hilendarski University of Plovdiv\\
236 Bulgaria Blvd.\\
4003 Plovdiv, Bulgaria
\\
e-mail: dobrinka@uni-plovdiv.bg}

\end{document}